\documentclass{amsart}

\usepackage{amssymb}
\let\cal\mathcal
\def\Ascr{{\cal A}}

\def\Cscr{{\cal C}}

\def\Hom{\operatorname {Hom}}

\def\r{\rightarrow}

\newtheorem{lemma}{Lemma}

\newtheorem{theorem}[lemma]{Theorem}

\theoremstyle{remark}

\newtheorem{remark}[lemma]{Remark}

\newdimen\uboxsep \uboxsep=1ex
\def\uboxn#1{\vtop to 0pt{\hrule height 0pt depth 0pt\vskip\uboxsep
\hbox to 0pt{\hss #1\hss}\vss}}

\def\uboxs#1{\vbox to 0pt{\vss\hbox to 0pt{\hss #1\hss}
\vskip\uboxsep\hrule height 0pt depth 0pt}}

\author{M. Van den Bergh}
 \address{Departement WNI,  Limburgs Universitair Centrum,
  3590 Diepenbeek, Belgium.}
    \email{vdbergh@luc.ac.be}
    \thanks{The  author is a senior researcher at the FWO}
    \date{\today}
    \title{A remark on a theorem by Deligne}
     
\begin{document}
     \begin{abstract}
     We give a proof avoiding spectral sequences of Deligne's decomposition
     theorem for objects in a triangulated category admitting 
     a Lefschetz homomorphism.
     \end{abstract}
     \maketitle
Below $\Ascr$ is a triangulated category equipped with
a bounded t-structure. In addition $\Ascr$ will be equipped with an
auto-equivalence $A\mapsto A(1)$
compatible with the t-structure. 

In \cite{de3,de2} Deligne proves the following result:
\begin{theorem} \label{del}Let $A$ be an object of $\Ascr$ equipped with
an endomorphism $\phi:A[-1]\r A[1](1)$ such that its iterates
$\phi^n:A[-n]\r A[n](n)$ induce isomorphisms $H^{-n}(A)\r H^{n}(A)(n)$.
Then there exists an isomorphism $A\cong \oplus_k H^{-k}(A)[k]$.
\end{theorem}
Deligne's proof is slightly indirect. He first shows that a decomposition
as asserted in Theorem \ref{del} exists if and only if for every cohomological 
functor
$T:\Ascr\r \Cscr$ to an abelian
category  the resulting spectral sequence 
\begin{equation}
\label{deg}
E_2^{pq}:T^p(H^q(A))\Rightarrow T^{p+q}(A)
\end{equation}
degenerates. He then proceeds to show that \eqref{deg} does indeed degenerate.

The aim of this note is to give a  proof of Theorem \ref{del}
which avoids the use of spectral sequences. 

Let $A,\phi$  be as in the statement of Theorem \ref{del}. We start
with the following statement: 
\begin{itemize}
\item[($\mathrm{Hyp}_n$)] $A\cong A_n\oplus \left( \oplus_{|k|>n}
H^{-k}(A)[k]\right)$ with $A_n\in \Ascr^{[-n,n]}$.  
\end{itemize} 
By
the boundedness of the t-structure ($\mathrm{Hyp}_n$) is true for $n\gg 0$. 
We need to show that it is true for $n=0$, so we use descending induction on
$n$.

Assume  ($\mathrm{Hyp}_n$) is true for $n\ge 1$. We will show that 
($\mathrm{Hyp}_{n-1}$) is also true. Without loss of generality we may 
assume that the isomorphism in the statement of ($\mathrm{Hyp}_n$)
is an equality. Let $i:A_n \r A$, $p:A\r A_n$ be respectively the inclusion
and the projection map. They induce identifications $H^q(A)=H^q(A_n)$ for
$|q|\le n$. Let $\alpha$ be the composition of the following maps
\[
H^{-n}(A)\r A_n[-n] \xrightarrow{i} A[-n]\xrightarrow{\phi^n} A[n](n)
\xrightarrow{p} A_n[n](n) \r H^{n}(A)(n)
\]
where the first and last map are obtained from the canonical maps 
$H^{-n}(A_n)[n]\r A_n \r H^{n}(A_n)[-n]$ which exist because 
$A_n\in \Ascr^{[-n,n]}$.

Applying $H^0(-)$ and the hypotheses we find that $\alpha=H^0(\phi^n)$
and hence it is an isomorphism. Composing  arrows we find that $\alpha$ is also
the  composition of maps
\begin{gather}
H^{-n}(A)\r A_n[-n] \r H^{n}(A)(n) \label{comp1}\\
H^{-n}(A)\r A_n[n] \r H^{n}(A)(n) \label{comp2}
\end{gather}
inducing isomorphisms on $H^0$.

From \eqref{comp1} it follows that 
\begin{equation}
\label{firstsplit}
A_n[-n]\cong H^{-n}(A)\oplus C
\end{equation}
 for
some $C\in \Ascr^{[1,2n]}$. 
Shifting and substituting in \eqref{comp2} we deduce that 
$\alpha$ is a composition 
\[
H^{-n}(A)\r 
H^{-n}(A)[2n]\oplus C[2n]
 \r H^{n}(A)(n) 
\]
Since $\Hom_\Ascr(H^{-n}(A)[2n], H^{n}(A)(n))=0$ we see that $\alpha$ is
actually a composition 
\[
H^{-n}(A)\r C[2n]\r H^{n}(A)(n)
\]
and these maps still induce isomorphisms in degree zero.
Thus $
C[2n]\cong H^{n}(A)(n)\oplus D
$ for $D\in \Ascr^{[-2n+1,-1]}$. Shifting and substituting in 
\eqref{firstsplit} yields a decomposition
\[
A_n\cong H^{-n}(A)[n] \oplus H^{n}(A)(n)[-n] \oplus D[-n] 
\]
Putting $A_{n-1}=D[-n]$ finishes the induction step and the proof.
\begin{remark} It follows from the above proof that the decomposition
asserted in Theorem \ref{del} still exists if we have maps
$\phi_n: A[-n]\r A[n](n)$ inducing isomorphisms $H^{-n}(A)\r H^{n}(A)(n)$ 
which are not necessarily powers of a fixed $\phi:A[-1]\r A[1](1)$. However
I have no example where this extra generality applies.
\end{remark}
\begin{remark}
In \cite{de2} Deligne constructs several canonical isomorphisms 
$A\cong \oplus_k H^{-k}(A)[k]$. We have not tried to duplicate
these constructions with our approach.
\end{remark}

\begin{thebibliography}{1}

\bibitem{de3}
P.~Deligne, {\em Th\'eor\`eme de {L}efschetz et crit\`eres de
  d\'eg\'en\'erescence de suites spectrales}, Inst. Hautes \'Etudes Sci. Publ.
  Math. (1968), no.~35, 259--278.

\bibitem{de2}
P.~Deligne, {\em D\'ecompositions dans la cat\'egorie d\'eriv\'ee}, Motives
  (Seattle, WA, 1991) (Providence, RI), Proc. Sympos. Pure Math., vol.~55,
  Amer. Math. Soc., Providence, RI, 1994, pp.~115--128.

\end{thebibliography}

\def\cprime{$'$}
\ifx\undefined\bysame
\newcommand{\bysame}{\leavevmode\hbox to3em{\hrulefill}\,}
\fi

\end{document}